\documentclass[12pt,a4paper]{amsart}

\usepackage {amsfonts}
\usepackage[english]{babel}
\def\g{\gamma}
\def\d{\delta}
\def\a{\alpha}
\def\b{\beta}
\def\e{\varepsilon}
\def\l{\lambda}
\def\R{{\mathbb R}}
\def\C{{\mathbb C}}
\def\N{{\mathbb N}}
\def\Z{{\mathbb Z}}
\def\Re{{\rm Re}}
\def\Im{{\rm Im}}
\def\dist{\mbox{dist}}

\def\supp{\mbox{supp}}

\newtheorem{theo}{Theorem}
\newtheorem{de}{Definition}

\newtheorem{Le}{Lemma}
\begin{document}

\title{Periodicity of almost periodic zero sets \\ and almost periodic functions}

\author{S.Yu. Favorov}

\address{Sergey Yuryevich Favorov,
\newline\hphantom{iii} The Karazin Kharkiv National University
\newline\hphantom{iii} Svoboda squ., 4,
\newline\hphantom{iii} 61022, Kharkiv, Ukraine}
\email{sfavorov@gmail.com}

\maketitle {\small
\begin{quote}
\noindent{\bf Abstract.} Whenever all differences between zeros of
two holomorphic almost periodic functions in a strip form a
discrete set, then both functions are infinite products of
periodic functions with commensurable periods. In particular, the
result is valid for some classes of Dirichlet series.
\medskip


\end{quote}
}

\bigskip
It was proved in \cite{FG} that each quasipolynomial
 \begin{equation}\label{Q}
   Q(z)=\sum_{n=1}^Na_ne^{\l_nz},\quad \l_n\in\R, \quad
   a_n\in\C.
 \end{equation}
with a discrete set of differences between its zeros is periodic
up to a multiplier without zeros, i.e., is of the form
 \begin{equation}\label{Q1}
Q(z)=Ce^{\b z}\prod_{k=1}^{N}\cosh(\omega
  z+b_k),\quad \b,\omega\in\R,
   \quad \ \C,\,
   b_k\in\C.
\end{equation}
The result is also valid for infinite sums
 \begin{equation}\label{S}
 S(z)=\sum_{n=1}^\infty a_ne^{\l_nz},\quad \l_n\in\R, \quad  a_n\in\C,
 \end{equation}
under conditions
 \begin{equation}\label{S1}
 \sum_{n=1}^\infty |a_n|<\infty, \quad \l_1=\sup_n\l_n<\infty,\quad \l_2=\inf_n\l_n>-\infty,\quad a_1a_2\neq0.
 \end{equation}
Note that zeros of  functions (\ref{Q}), (\ref{S}) are located
in a vertical strip of a finite width.

In 1949 M.G. Krein and B.Ja.Levin (see \cite{KL}, also \cite{L},
 ch.6, p.2 and Appendix~6) introduced and studied the class
$\Delta$ of entire almost periodic functions of an exponential
growth with zeros in a horizontal strip of a finite width. In
particular, if a sum $S(z)$ in (\ref{S}) satisfies condition
(\ref{S1}), then $S(iz)$ belongs to $\Delta$. In fact,
representation (\ref{Q1}) was obtained in \cite{FG} for functions
from $\Delta$ with a discrete sets of differences between zeros
(of course, it should be replaced $e^{\b z}$ by $e^{i\b z}$ and
$\cosh(\omega z+b_k)$ by $\cos(\omega z+b_k)$ ).

The phenomenon takes place for a pair of functions and discreteness of the set of differences between their zeros as well
(\cite{G}). For example, let $Q_1,\,Q_2$ be entire functions of the form (\ref{S}) under conditions (\ref{S1}).
 If the set $\{z-w,\ Q_1(z)=0,\,Q_2(w)=0\}$ is discrete, then the both functions have
the form (\ref{Q1}) with the same $\omega$ and possibly different $C,\b,N,b_k$. In particular, in the case $Q_2(-z)=Q_1(z)=Q(z)$ we
get representation (\ref{Q1}) for any function (\ref{S}) with a discrete set $\{z+z',\ Q(z)=Q(z')=0\}$.

The proof of the above results are based on the  property of zeros
$\{z_j\}$ of functions $f\in\Delta$ (\cite{L}, Appendix~6, p.2) to
be almost periodic in the  sense of the following definition:
\begin{de}
A zero set $\{z_j\}$ is almost periodic if for any $\e>0$ there is a
relatively dense set $E_\e\subset\R$  such that for each $\tau\in
E_\e$ there exists a bijection $\sigma:\,\N\to\N$ with the
property
 \begin{equation}\label{ap2}
\sup_j |z_j+i\tau-z_{\sigma(j)}|<\e.
 \end{equation}
 \end{de}
 Recall that a set
$E\subset\R$ is relatively dense, if there exists $L<\infty$ such
that $E\cap[a,a+L]\neq\emptyset$ for any $a\in\R$.

\medskip
Practically, the above definition is applied to zero sets in a {\it closed} vertical strip.
In the present paper we investigate holomorphic
almost periodic functions and almost periodic sets in an open
strip, in particular, in the complex plane and in half-planes.

\begin{de}
A continuous function $f(z)$ in a strip $S=\{z:\,a<\Re z<b\},\
-\infty\le a<b\le\infty$ is almost periodic, if for any substrip
$S_0$, $\overline{S_0}\subset S$ and any $\e>0$ there is a
relatively dense set $E_{\e,S^0}\subset\R$ such that for any
$\tau\in E_{\e,S^0}$
 \begin{equation}\label{apf}
 \sup_{z\in S_0}|f(z+i\tau)-f(z)|<\e.
 \end{equation}
\end{de}
A typical example of an entire almost periodic function is sum
(\ref{S}) under  conditions
 $$
 \sum_{n=1}^\infty |a_n|<\infty,\quad \sup_n\l_n<\infty,\quad \inf_n\l_n>-\infty.
 $$
  Moreover, a holomorphic function $f(z)$ in a strip
$S=\{z:\,a<\Re z<b\},\ -\infty\le a<b\le\infty$ is almost periodic
if and only if there exists a sequence of quasipolynomials $Q_n$
of the form (\ref{Q}) such that  for any substrip $S_0$,
$\overline{S_0}\subset S$,
 $$
 \sup_{z\in S_0}|f(z)-Q_n(z)|\to0,\quad n\to\infty.
 $$
(see \cite{JT}, item 8).

In order to take into account multiplicities of zeros, we  use
the term {\it divisor} instead of zero set. Namely, a divisor $Z$ in
a domain $D$ is a mapping $Z:\,D\to\N\cup\{0\}$ such that
$|Z|=\supp Z$ is a set without limit points in $D$. In other
words, a divisor in $D$ is the sequence of points $\{z_j\}\subset
D$ that has no limit points in $D$ and every point of $D$ appears at
most a finite times in the sequence. If $Z(z)\le1$ for all $z\in D$,
we will identify $Z$ and $|Z|$. The divisor of zeros of a
holomorphic function $f$ in $D$ is the map $Z_f$ such that
$Z_f(a)$ equals multiplicity of zero of function $f$ at the point
$a$.

The following definition appeared at first in \cite{T}:

\begin{de}[\cite{T},\cite{FRR}]\label{d3} A divisor $Z=\{z_j\}$ in a strip
$S=\{z:\,a<\Re z<b\},\ -\infty\le a<b\le\infty$, is called almost
periodic if for any $\e>0$ and any substrip $S^0$,
$\overline{S^0}\subset S$, there is a relatively
dense set
 \begin{equation}\label{ap3}
 E_{\e,S^0}=\{\tau\in\R:\,z_j\in S^0 \bigvee z_{\sigma(j)}\in S^0
\Longrightarrow |z_j+i\tau-z_{\sigma(j)}|<\e\}.
 \end{equation}
where $\sigma=\sigma_\tau$ is a suitable bijection $\N\to\N$.
\end{de}

{\bf Remark.} If we take $\sigma^{-1}(j)$ instead of $j$ in
(\ref{ap3}), we get
\begin{equation}\label{ap1}
z_j\in S^0 \bigvee z_{\sigma^{-1}(j)}\in S^0 \Longrightarrow
|z_j-i\tau-z_{\sigma^{-1}(j)}|<\e.
 \end{equation}

\begin{theo}[\cite{FRR}]\label{t1}
a) If $f$ is a holomorphic function in a strip $S$ with an almost
periodic modulus, then the divisor of zeros of $f$ is almost
periodic,

b) for any almost periodic divisor $Z$ in a strip $S$ there is a
holomorphic function $f$ in $S$ with an almost periodic modulus
such that $Z=Z_f$,

c) for any almost periodic divisor $Z$ in a strip $S$ such that
$|Z|\cap S'=\emptyset$ for some open substrip $S'\subset S$ there
is a holomorphic almost periodic function $f$ in $S$ such that
$Z=Z_f$.
\end{theo}

 \medskip
 There exist almost periodic divisors in the plane with a
 discrete set of differences that is not periodic.

 {\bf Example 1.} Let
 $$
 Z=\{z_{n,k}=i2^{nk}+2^k,\ n\in\Z,\ k\in\N\}
 $$
be a discrete set in $\C$. It is easy to see that the set is
almost periodic, differences between its points form a discrete
set, but $Z$ is a countable union of periodic sets with different
commensurable periods.

\medskip
In our article we prove the following theorem

\begin{theo}\label{t2}
Let $Z$, $W$ be  almost periodic divisors in a strip
$S=\{z:\,a<\Re z<b\},\ -\infty\le a<b\le\infty$. If

a) $|Z|\cap S'=\emptyset$ for some substrip $S'\subset S$,

b)  for any substrip $S_0$ of a finite width,
$\overline{S_0}\subset S$, the set $\{z-w,\ z\in |Z|\cap
S_0,\,w\in |W|\cap S_0\}$ is discrete,

\noindent then $Z,\,W$ are at most countable sums of periodic
divisors with commensurable periods. \end{theo}

Show that condition a) is essential.

{\bf Example 2.} Let $S=\{z:\,|\Re z|<1\},\
Z=\{z_{m,n}=(m+in)e^{i\a}\in S,\,m,n\in\Z\}$, where $\a\in\R$ such
that $\cot\,\a$ is an irrational number. Clearly, the set of
differences of elements of $Z$ is discrete. Let us prove that $Z$
is an almost periodic divisor without any periods.

By Kronecker Lemma (see, for example, \cite{Lt}, Ch.2, \$ 2), for
any $\e>0$ the inequalities
$$
|t\cot\a|<\e\ (\text{mod}\,\Z),\quad |t|<\e\ (\text{mod}\,\Z)
 $$
 has a relatively dense set of
common solutions. Therefore, the inequality
$$
|m\cot\a|<\e(1+|\cot\a|)\ (\text{mod}\,\Z)
$$
has a relatively dense set of integer solutions. The later means
that for any $\d$ there exist pairs of integers $(m,n)\in\Z^2$
such that
 \begin{equation}\label{in}
\quad |m\cos\a-n\sin\a|<\d,
\end{equation}
and the set of $m$ with this property is relatively dense. Put
 $$
E=\{m\sin\a+n\cos\a:\,|m\cos\a-n\sin\a|<\d\}.
 $$
Let $\tau=\tilde m\sin\a+\tilde n\cos\a\in E$. For any $z_{m,n}\in
Z$ and $m'=m+\tilde m,\ n'=n+\tilde n$ we have
 $$
z_{m,n}+i\tau-z_{m',n'}=(m+in)e^{i\a}+i\Im[(\tilde m+i\tilde
n)e^{i\a}]-(m'+in')e^{i\a}=\Re[(\tilde m+i\tilde n)e^{i\a}].
 $$
By (\ref{in}), $|z_{m,n}+i\tau-z_{m',n'}|<\d$. Since
$$
|\tau-(\cos^2\a/\sin\a+\sin\a)\tilde m|=|\tilde n\cos\a-(\cos^2\a/\sin\a)\tilde m| <\d|\cot\a|,
$$
 we see that the set $E$ is
relatively dense. So, $Z$ is an almost periodic divisor.

 Since $Z$ sites in the vertical strip of width $2$, we see
that any period of $Z$ must have the form $iT,\ T\in\R$. Hence,
for some $(n,m),\,(n',m')\in\Z^2,\ (n,m)\neq(n',m')$,
 $$
(n+im)e^{i\a}-(n'+im')e^{i\a}=iT,
 $$
and
$$
(n-n')+i(m-m')=iT(\cos\,\a-i\sin\,\a).
 $$
This equality contradicts our choice of $\a$. Hence any part of
$Z$ has no periods.

\medskip
Our proof of Theorem \ref{t2} makes into use the following Lemmas

\begin{Le}\label{L1}
Let $Z=\{z_n\}$ be an almost periodic divisor in a strip $S$,
$S^0$ be a substrip, $\overline{S_0}\subset S,\ \e>0,\
\tau_1,\,\tau_2\in E_{\e,S^0}$, where $E_{\e,S^0}$ is from Definition \ref{d3}.
 Then there is a bijection $\sigma:\N\to\N$ such that
if $z_j\in S^0,\,\dist(z_j,\partial S^0)>\e$, then
\begin{equation}\label{ap0}
 |z_j+i(\tau_1-\tau_2)-z_{\sigma(j)}|<2\e,\quad
 |z_j-i(\tau_1-\tau_2)-z_{\sigma^{-1}(j)}|<2\e.
\end{equation}
 \end{Le}
It can be proved easily that the similar assertion is valid for
$\tau_1+\tau_2$ as well.

{\bf Proof.} If $z_j$ satisfies conditions of the Lemma, then for
some bijection $\sigma_1:\N\to\N$ we have
$|z_j+i\tau_1-z_{\sigma_1(j)}|<\e$, therefore
$|\Re(z_{\sigma_1(j)}-z_j)|<\e$ and $z_{\sigma_1(j)}\in S^0$. By
(\ref{ap1}), we also get the inequality
$|z_{\sigma_2^{-1}\circ\sigma_1(j)}-i\tau_2-z_{\sigma_1(j)}|<\e$ for some bijection $\sigma_2:\N\to\N$.
Therefore, we obtain the first part in (\ref{ap0}) with
$\sigma=\sigma_2^{-1}\circ\sigma_1$. If we change places of $\tau_1$ and
$\tau_2$, we obtain the second part in (\ref{ap0}).

 \begin{Le}\label{L2} Let $Z=\{z_n\},\,W=\{w_m\}$ be almost periodic
divisors in a strip $S$.  Then for any $\e>0$ and substrip
$S^0\subset S$, $\dist(S_0,\partial S)>\e$, there is a relatively
dense set $E\subset\R$ with the  property:  for any $\tau\in E$
there are bijections $\sigma_Z,\,\sigma_W$ from $\N$ to $\N$ such
that for all $z_j\in |Z|\cap S^0,\,w_r\in |W|\cap S^0$
\begin{equation}\label{ap4}
 |z_j+i\tau-z_{\sigma_Z(j)}|<\e,\quad
 |w_r+i\tau-w_{\sigma_W(r)}|<\e
 \end{equation}
 and
\begin{equation}\label{ap5}
 |z_j-i\tau-z_{\sigma^{-1}_Z(j)}|<\e,\quad
 |w_r-i\tau-w_{\sigma^{-1}_W(r)}|<\e.
 \end{equation}
 \end{Le}

{\bf Proof.} Put $\tilde S=\{z\in S:\,\dist(z,S^0)<\e/2\}$. By
(\ref{ap3}) and (\ref{ap1}), there is a number $L$ such that any
interval of length $L$ contains $\tau_Z$ and $\tau_W$ with the
properties
 \begin{equation}\label{ap6}
 |z_j+i\tau_Z-z_{\sigma_Z(j)}|<\e/4,\quad
 |z_j-i\tau_Z-z_{\sigma^{-1}_Z(j)}|<\e/4
 \end{equation}
 for all $z_j\in |Z|\cap\tilde S$, and
 \begin{equation}\label{ap7}
 |w_r+i\tau_W-w_{\sigma_W(r)}|<\e/4,\quad \quad
 |w_r-i\tau_W-w_{\sigma^{-1}_W(r)}|<\e/4
 \end{equation}
 for all $w_r\in |W|\cap\tilde S$. Here $\sigma_Z,\,\sigma_W$ are some
 (in general, different) bijections $\N\to\N$.

 We may suppose that $N=2L/\e$ is an integer. Hence for
any $k\in\Z$ there are integers $n(k),\,m(k),\,0\le n(k),m(k)\le
N$ such that
 $$
 |kL+n(k)\e/2-\tau_Z|<\e/4,\quad |kL+m(k)\e/2-\tau_W|<\e/4.
 $$
  The differences $n(k)-m(k)$ take at most $2N+1$ values. Choose
$k_1,\dots,k_r,\,r\le 2N+1,$ such that for any $k\in\Z$ there is
$k_s$ such that
 $$
 n(k)-m(k)=n(k_s)-m(k_s).
 $$
It follows easily that any interval of length $L(\max_s(|k_s|+2)$
contains a point of the form
 $$
 \tau=Lk+n(k)\e/2-Lk_s-n(k_s)\e/2=Lk+m(k)\e/2-Lk_s-m(k_s)\e/2.
 $$
If we replace $\tau_Z$ by $Lk+n(k)\e/2$ or $Lk_s+n(k_s)\e/2$ in
(\ref{ap6}), we obtain that this  inequality  satisfies with $\e/2$ instead
of $\e/4$.  The same is true if we replace $\tau_W$ by
$Lk+m(k)\e/2$ or $Lk_s+m(k_s)\e/2$ in (\ref{ap7}). Applying Lemma
\ref{L1} with $S^0=\tilde S$ and $\e/2$ instead of $\e$, we see
that $\tau$ satisfies (\ref{ap4}) and (\ref{ap5}).

 {\bf Proof of Theorem \ref{t2}.} Set a sequence of vertical substrips
$S_k,\,\overline{S_k}\subset S_{k+1},$ such that
 $$
|Z|\cap S_1\neq\emptyset,\quad |W|\cap S_1\neq\emptyset,\quad
S_1\cap S'\neq\emptyset,\quad\cup_k S_k=S.
 $$
 Take any points $z\in |Z|\cap S_1$, $w\in |W|\cap S_1$. It follows from
(\ref{ap4}) with $\e<\min\{\dist(z,\,\partial
S_1),\dist(w,\,\partial S_1)\}$ that there is $R<\infty$ such that
any horizontal strip of the width $R$ contains at least one point
$|Z|\cap S_1$ and at least one point $|W|\cap S_1$. By condition
b) of the theorem,  for each $k$ the set $\{z-w,\ z\in |Z|\cap
S_k,\,w\in |W|\cap S_k\}$ is discrete.  Hence there exists
$\g_k\in(0,\,1/2)$ such that whenever
$$
z_n,z_{n'}\in |Z|\cap S_{k+1},\ w_m,w_{m'}\in |Z|\cap S_{k+1},\quad z_n-w_m\neq z_{n'}-w_{m'},
$$
$$
|\Im(z_n-w_m)|<2R+3,\quad|\Im(z_{n'}-w_{m'})|<2R+3,
$$
 we get $\g_k<|(z_{n'}-w_{m'})-(z_n-w_m)|$. In particular, if we
put $w_m=w_{m'}$, then we get $\g_k<|z_n-z_{n'}|$ for any
$z_n,z_{n'}\in |Z|\cap S_{k+1}$, $z_n\neq z_{n'}$.

Fix $z_n\in |Z|\cap S_k$, and let a  number $\tau>1$ satisfies
(\ref{ap4}) for $Z$ and $W$ with $\e=\g_k/2$. Then there is a
unique
 $z_{n'}\in |Z|$ such that $|z_n+i\tau-z_{n'}|<\g_k/2$.
Indeed, otherwise we obtain
 $$
 |z_{n'}-z_{n''}|\le |z_n+i\tau-z_{n'}|+|z_n+i\tau-z_{n''}|<\g_k.
 $$
Set $T_k=(z_{n'}-z_n)/i$.

 Let $w_m\in |W|$ be such that
$|\Im\,(w_m-z_n)|<2R+2$. By (\ref{ap4}), there is a point
$w_{m'}\in |W|$ such that $|w_m+i\tau-w_{m'}|<\g_k/2$. Therefore,
  $$
|(z_n-w_m)-(z_{n'}-w_{m'})|\le|w_m+i\tau-w_{m'}|+|z_{n'}-z_n-i\tau|<\g_k.
  $$
   Since
  $$
|\Im(z_{n'}-w_{m'})|\le|\Im(z_n-w_m)|+|z_n-z_{n'}+i\tau|+|w_m-w_{m'}+i\tau|<2R+3,
 $$
  we get $z_n-w_m=z_{n'}-w_{m'}$ due to the choice of $\g_k$.
 Therefore, $w_{m'}=w_m+iT_k$.

 The latter equality takes place for all points
$|W|\cap\{w:\,\Im z_n-2R<\Im w<\Im z_n+2R\}$, in particular, for
some $w_l$ such that $\Im z_n+R<\Im w_l<\Im z_n+2R$. Namely, there
is $w_{l'}\in |W|$ such that $w_{l'}=w_l+iT_k$. Let $\zeta\in |Z|$
be any point from the set
\begin{equation}\label{set}
\{z:\,\Im z_n\le\Im z<\Im z_n+3R\}\subset\{z:\,\Im w_l-2R<\Im
z<\Im w_l+2R\}.
\end{equation}

 By (\ref{ap4}), there is a point
$\zeta'\in |Z|$ such that $|\zeta+i\tau-\zeta'|<\g_k/2$.
Therefore,
  $$
|(\zeta-w_l)-(\zeta'-w_{l'})|\le|\zeta+i\tau-\zeta'|+|iT_k-i\tau|<\g_k.
  $$
   Since $|\Im\zeta-\Im w_l|<2R$ and
  $$
|\Im(\zeta'-w_{l'})|\le|\Im(\zeta-w_l)|+|\zeta+i\tau-\zeta'|+|iT_k-i\tau|<2R+3,
 $$
  we get $\zeta-w_l=\zeta'-w_{l'}$ due to the choice of $\g_k$.
 Therefore, $\zeta'=\zeta+iT_k$.

In particular, there is a point $z_s\in |Z|\cap\{z:\,\Im
z_n+2R\le\Im z<\Im z_n+3R\}$ such that $z_s+iT_k\in |Z|$.
Continuing the line of reasoning, we obtain that for all $z\in
|Z|$ such that $\Im z\ge\Im z_n$ we get $z+iT_k\in |Z|$ and for
all $w\in |W|$ such that $\Im w\ge\Im z_n$ we get $w+iT_k\in |W|$.

If we take $w'_l\in |W|$ such that $\Im z_n-2R<\Im w'_l<\Im
z_n-R$, we also can find $w'_{l'}\in |W|$ such that
$w'_{l'}=w'_l+iT_k$. Next, we prove that for any point $\tilde
\zeta\in |Z|\cap\{z:\,\Im z_n-3R\le\Im z<\Im z_n\}$ there is
$\tilde\zeta'\in |Z|$ such that $\tilde\zeta'=\tilde\zeta+iT_k$.
Arguing as above, we show that for all $z\in |Z|$ such that $\Im
z\le\Im z_n$ we get $z+iT_k\in |Z|$ and for all $w\in |W|$ such
that $\Im w\le\Im z_n$ we get $w+iT_k\in |W|$.

Next, by (\ref{ap5}), take for any $z\in |Z|$ a point  $z''\in
|Z|$ such that $|z''+i\tau-z|<\g_k/2$. Then $z''+iT_k\in |Z|$ and
 $$
|(z''+iT_k)-z|\le |z''+i\tau-z|+|i\tau-iT_k|<\g_k.
 $$
Therefore, $z''+iT_k=z$ and $z-iT_k\in |Z|$ for all $z\in |Z|$. By
the same arguments, $w-iT_k\in |W|$ for all $w\in |W|$.

If an imaginary part of $T_k$ does not vanish, then  either
$z_n+iMT_k\in S'$, or $z_n-iMT_k\in S'$ for a suitable $M\in\N$. Therefore, $T_k$ is real and for any $M\in\Z$ we get
$(|Z|\cap S_k)+iMT_k=|Z|\cap S_k,\ (|W|\cap S_k)+iMT_k=|W|\cap
S_k$. Hence, the restrictions $Z\mid_{S_k},\,W\mid_{S_k}$ the
divisors $Z,\,W$ to $S_k$ are periodic divisors with period
$iT_k$.

The same arguments work for every $k=1,2,\dots$.

Let $T^0_k$ be the minimal common period of $Z\mid_{S_k}$ and
$W\mid_{S_k}$. Clearly, $T_k/T^0_k\in\N$. Besides, since $|Z|\cap
S_k\subset |Z|\cap S_m$ for $m>k$, we have $T_m/T^0_k\in\N$ as
well.

Finally, let $Z_k=Z\mid_{S_{k+1}\setminus S_k},\
W_k=W\mid_{S_{k+1}\setminus S_k}$. Then we obtain
 $$
 Z=Z_1+Z_2+Z_3+\dots,\qquad W=W_1+W_2+W_3+\dots.
 $$
Theorem is proved.

Theorem \ref{t2} implies the corresponding result for almost
periodic holomorphic functions.
\begin{theo}\label{t3}
Let $f$, $g$ be  almost periodic functions in a strip
$S=\{z:\,a<\Re z<b\}, -\infty\le a<b\le\infty$. If

a) either $f$, or $g$ has no zeros in an open substrip $S'\subset
S$,

b)  for any substrip $S_0$, $\overline{S_0}\subset S$, the set
$\{z-w,\ z,w\in S_0,\,f(z)=g(w)=0\}$ is discrete, then
\begin{equation}\label{pro}
f(z)=f_0(z)\prod_{k=1}^\infty f_k(z),\quad
g(z)=g_0(z)\prod_{k=1}^\infty g_k(z)
\end{equation}
where  $f_0,\,g_0$ are holomorphic almost periodic functions in
$S$ without zeros,  $f_k,\,g_k,\ k=1,2,\dots,$
 are periodic holomorphic in $S$ with commensurable periods $T_k$.
\end{theo}

{\bf Proof.} By Theorem \ref{t1}a), the divisors $Z,\,W$ of zeros of $f,\,g$, respectively, 
 are almost periodic and satisfy other conditions of the previous theorem. Let
$S_k,\,Z_k,\,W_k$, be the same as
in the proof of Theorem \ref{t2}. Set
 $$
 S_k=\{z:\,\eta_k<\Re\,z<\eta'_k\},\quad\eta_k<\eta_{k-1},\quad\eta'_k>\eta'_{k-1},\ \forall k.
 $$
For any $k$ there is only a finite number of points
$a^k_1,\dots,a^k_{m_k}\in |Z_k|\cup\{z:\,0\le\Im z
<T_k\}$.  Therefore, $|Z_k|=\{a_1^k,\dots,a_{m_k}^k\}+iT_k\Z$. Take
$\e_k^j>0,\,1\le j\le m_k,\,3\le k<\infty,$ such that
 \begin{equation}\label{s}
\sum_{j,k}\e_k^j<\infty.
 \end{equation}
Fix $k>1$. Let
 $$
 h_j(w)=1-w\exp(-2\pi a^k_j /T_k ),\quad
j=1,\dots,m_k.
$$
Clearly, $h_j(\exp(2\pi z/T_k))$ has a divisor $a_j^k+iT_k\Z$.
Since $\a_j^k\in S_k\setminus S_{k-1`}$, we see that either $\eta'_{k-1}\le\Re\,a_j^k$, or
$\eta_{k-1}\ge\Re\,a_j^k$. In the first case,  $\log h_j(w)$ is
holomorphic on the disc $|w|<\exp(2\pi\Re\,a^k_j /T_k)$, and there
is a polynomial $P_j^k$  such that
 $$
|\log h_j(w)-P_j^k(w)|<\e_k^j\quad\text{for}\  |w|\le\exp(2\pi
\eta'_{k-2}/T_k).
 $$
In the second one, $\log h_j(w)$ is holomorphic on the set
$|w|>\exp(2\pi\Re\,a^k_j /T_k)$, and there is a polynomial $P_j^k$
 such that
 $$
|\log h_j(w)-P_j^k(w^{-1})|<\e_k^j\quad\text{for}\
|w|\ge\exp(2\pi \eta_{k-2}/T_k).
 $$
Put $Q_j(z)=P_j^k(\exp(2\pi z/T_k))$ in the first case, and
$Q_j(z)=P_j^k(\exp(-2\pi z/T_k))$ in the second one. We obtain
 \begin{equation}\label{apr}
|h_j(e^{2\pi z/T_k})e^{-Q_j(z)}|<e^{\e_j^k}\quad\text{for}\  z\in
S_{k-2}.
 \end{equation}
Put
 $$
 f_k(z)=\prod_{j=1}^{m_k}h_j(e^{2\pi z/T_k}),\ k=1,2,\quad
 f_k(z)=\prod_{j=1}^{m_k}h_j(e^{2\pi z/T_k})e^{-Q_j(z)},\ k>2.
 $$
Clearly, $Z_{f_k}=Z_k$. By (\ref{s}) and (\ref{apr}), the product
in (\ref{pro}) converges in every substrip $S_k$. The function
$f(z)/\prod_{k=1}^\infty f_k(z)$ is holomorphic in $S$, and, by
\cite{T}, Theorem 1, it is almost periodic in $S$. In the same way,
we obtain the representation of $g$.

{\bf Remark.} It follows easily that for entire almost periodic
functions $f,\,g$ with zeros in a strip $\tilde S$ of a finite
width one can take $S_1=\tilde S$. Therefore, $f,\,g$ are periodic
functions with the same period up to almost periodic multiplies
without zeros. Hence the results of \cite{FG}, \cite{G} follow
from Theorem \ref{t3}.

\end{document}